\documentclass[12pt]{article}
\usepackage{amssymb}
\usepackage{amsmath}
\usepackage{indentfirst,latexsym}
\usepackage{mathrsfs}
\usepackage{graphicx}
\usepackage{fancyhdr}
\date{}
\begin{document}

\title{The descent statistic on involutions is not log-concave}
\author{Marilena Barnabei, Flavio Bonetti, and Matteo Silimbani \thanks{
Dipartimento di Matematica - Universit\`a di Bologna}} \maketitle

\noindent {\bf Abstract.} We establish a combinatorial connection
between the sequence $(i_{n,k})$ counting the involutions on $n$
letters with $k$ descents and the sequence $(a_{n,k})$ enumerating
the semistandard Young tableaux on $n$ cells with $k$ symbols.
This allows us to show that the sequences $(i_{n,k})$ are not
log-concave for some values of $n$, hence answering a conjecture
due to F. Brenti.
\newline

\noindent {\bf Keywords:} involution, descent, semistandard Young
tableau, reverse Yamanouchi word.

\noindent {\bf AMS classification:} 05A05, 05A15, 05A19, 05E10.

\section{Introduction}

\noindent A finite sequence $(s_1,\ldots,s_n)$ of real numbers is
said to be \emph{unimodal} if there exists an index $t$ such that
$s_1\leq s_2\leq\cdots\leq s_t$ and $s_t\geq\cdots \geq
s_{n-1}\geq s_n$. An arbitrary sequence $(s_i)_{i\in\mathbb{N}}$
is \emph{log-concave} if $s_{i-1}\cdot s_{i+1}\leq s_i^2$ for
every $i>0$. It is immediately seen that a finite log-concave
sequence of positive numbers is unimodal. In recent years, several
authors focused on the study of such two properties in relation to
the distribution of the descent statistic on involutions. More
precisely, given a word $w=w_1\ldots w_n$ on a linearly ordered
alphabet, the \emph{descent set} of $w$ is defined as
des$(w)=\{1\leq i<n:w_i\geq w_{i+1}\}$ and  the cardinality of the
set des$(w)$ is denoted by $d(w)$. An analogous definition can be
given for the \emph{ascent set} of a word. If $\sigma$ is an
involution, the descent set of $\sigma$ is the descent set of the
word $\sigma(1)\ldots\sigma(n)$. Let $i_{n,k}$ be the number of
involutions on $n$ letters with $k$ descents and let
$$I_n(x)=\sum_{k=0}^{n-1} i_{n,k}x^k$$ be the generating function of
the sequence $i_{n,k}$, for every $n\in\mathbb{N}$. Strehl
\cite{stre} proved that the coefficients of $I_n(x)$ are
symmetric. Recently, Brenti (see \cite{duk}) conjectured that the
coefficients of the polynomial $I_n(x)$ are log-concave. Dukes
\cite{duk} obtained some partial results on the unimodality of
such coefficients and Guo and Zeng \cite{zeng} succeeded in
proving that the sequence $i_{n,k}$ is unimodal.

\noindent In the present note, we disprove Brenti's conjecture
exploiting the combinatorial relation of the sequence $i_{n,k}$
with the sequence $a_{n,s}$ counting semistandard Young tableaux
on $n$ cells with $s$ symbols. This last sequence is easily seen
to be not log-concave. We show that the generating functions of
these two sequences are related by a binomial transformation. This
fact allows us to refute the log-concavity of the polynomial
$I_n(x)$. Moreover, we deduce an explicit formula for the integers
$i_{n,k}$.\newline

\noindent The relation between standard and semistandard Young
tableaux sheds new light on the combinatorial properties of Young
tableaux. For example, it provides an immediate proof of the well
known fact that every Schur function $s_{\lambda}$ can be
expressed as a sum of suitable fundamental quasi-symmetric
functions. Moreover, the present techniques allow to investigate
other properties of the distribution of the descent statistic both
on the set of involutions itself \cite{bbsb} and on some notable
subset of involutions \cite{bbsa}.

\section{Standard Young Tableaux}

\noindent Consider the set $\mathscr{T}_n$ of standard Young
tableaux on $n$ cells. It is well known that the
Robinson-Schensted algorithm establishes a bijection
$\psi:\mathscr{I}_n\to \mathscr{T}_n$, where $\mathscr{I}_n$ is
the set of involutions over $[n]:=\{1,2,\ldots,n\}$. A further
bijective map $\chi$ exists between $\mathscr{T}_n$ and the set
$\mathscr{Y}_n$ of reverse Yamanouchi words of length $n$. We
recall that a \emph{reverse Yamanouchi word} is a word $w$ with
integer entries such that any left subword of $w$ does not contain
more occurrences of the symbol $(i+1)$ than of $i$, for every
$i\geq 1$. The Yamanouchi word $\chi(T)$ associated to a given
tableau $T$ is obtained by placing in the $i$-th position the row
index of the cell of $T$ containing the symbol $i$. For example,
if $T$ is the standard Young tableau

$$T=\begin{array}{cccc}
1&3&5& \\2&6&7&\\4&8&&
\end{array}$$
we have $\chi(T)=1\,2\,1\,3\,1\,2\,2\,3$.\newline

\noindent Clearly, the composition $\varphi:=\chi\circ\psi$ yields
a bijection between the sets $\mathscr{I}_n$ and $\mathscr{Y}_n$.
We remark that $\varphi$ turns each ascent of a given involution
$\sigma$ into a descent of the correspondent reverse Yamanouchi
word.\newline

\noindent Let $i_{n,h}$ be the number of involutions $\sigma\in
\mathscr{I}_n$ with $h$ descents and $y_{n,k}$ the number of
reverse Yamanouchi words of length $n$ with $k$ descents. The
preceding remark implies that
$$y_{n,k}=i_{n,n-1-k}.$$
The present approach leads to an immediate proof of the following
result originally due to Strehl \cite{stre}:
\newtheorem{yama}{Proposition}
\begin{yama}
\label{lasimm} For every $n\in\mathbb{N}$, we have
$$y_{n,k}=y_{n,n-1-k}.$$
\end{yama}

\noindent \emph{Proof} Given a reverse Yamanouchi word $y$,
consider the conjugate word $\tilde{y}$ defined as follows: if
$y_i=m$, then $\tilde{y}_i$ is the number of occurrences of the
integer $m$ in the left subword $y_1\ldots y_i$. For example, the
conjugate of the word $$y=1\,2\,1\,3\,1\,2\,2\,3$$ is
$$\tilde{y}=1\,1\,2\,1\,3\,2\,3\,2$$
Note that, if $y$ is associated with the tableau $T$, $\tilde{y}$
is associated with the conjugate tableau of $T$. Clearly, $y$ has
$k$ descents if and only if $\tilde{y}$ has $n-1-k$
descents.\begin{flushright} $\diamond$
\end{flushright}

\noindent For every $n\in\mathbb{N}$, define
$$I_n(x)=\sum_{\sigma\in \mathscr{I}_n} x^{d(\sigma)}=\sum_{h=0}^{n-1} i_{n,h}x^h.$$
This polynomial can be rewritten in terms of reverse Yamanouchi
words as follows:
$$I_n(x)=\sum_{k=0}^{n-1} y_{n,n-1-k}x^k=\sum_{k=0}^{n-1} y_{n,k}x^k=\sum_{y\in \mathscr{Y}_n} x^{d(y)}.$$

\section{Semistandard Young Tableaux}
\label{serve}

\noindent Given a Ferrers diagram $\lambda$, a \emph{semistandard
Young tableau} on $k$ symbols of shape $\lambda$ is an array
obtained by placing into each cell of the diagram an integer in
$[k]$ so that the entries are strictly increasing by rows and
weakly increasing by columns. We consider the infinite matrix
$A=(a_{n,k})$, with $n,k\in\mathbb{N}$, where $a_{n,k}$ denotes
the number of semistandard Young tableaux with $n$ cells and $k$
symbols. An explicit expression for the column generating function
$F_k(x)$ of the matrix $A$
$$F_k(x)=\sum_{n\geq 0}a_{n,k}x^n=\frac{1}{(1-x)^k(1-x^2)^{{k\choose 2}}}$$
was firstly given by Schur (see \cite{Lot}). This yields
immediately the following explicit formula for the integers
$a_{n,k}$:

\begin{equation}a_{n,k}=\sum_{j=0}^{\lfloor \frac{n}{2}\rfloor}{\frac{k(k-1)}{2}+j-1
\choose j}{k+n-2j-1\choose k-1}.\label{a}\end{equation}

\noindent  The following properties of the row sequences of the
matrix $A$ are direct consequences:

\newtheorem{eazy}[yama]{Proposition}
\begin{eazy}
\label{lapide} The sequence $(a_{n,k})_{k\in\mathbb{N}}$ is in
general not log-concave.
\end{eazy}

\noindent \emph{Proof} Exploiting Formula (\ref{a}), we have:
$$a_{45,2}^2=304704<307970=a_{45,1}\cdot a_{45,3}.$$
\begin{flushright} $\diamond$
\end{flushright}

\noindent We are now interested in establishing a connection
between the sequences $(a_{n,k})$ and $(y_{n,k})$. To this aim, we
associate with a given semistandard tableau $T$ a biword $(w,y)$
as follows: $w$ contains all the entries in $T$ listed in
non-decreasing order. The word $y$ is obtained by listing the row
indices of the occurrences of each symbol, starting from the
smallest one. If a symbol $j$ occurs more than once, we write the
corresponding row indices in increasing order. It is easy to check
that $y$ is a Yamanounchi word. Note that the biword $(w,y)$
uniquely determines the tableau $T$. In fact, applying to the
biword $(w,y)$ the Robinson-Schensted-Knuth column insertion
procedure, we get the pair $(Y,T)$, where $Y$ is a row Yamanouchi
tableau, namely, a tableau whose $i$-th row consists only of
letters $i$ for all $i$.
\newline

\noindent We are now going to show that the number of semistandard
tableaux on $s$ symbols associated with a given reverse Yamanouchi
word $y$ depends only on the number of descents of $y$. Fix a
reverse Yamanouchi word $y$ with $k$ descents. Any semistandard
tableau with associated biword $(w,y)$ must contain at least $k+1$
different symbols. In fact, if $y$ has a descent at position $i$,
by the definition of the correspondence between tableaux and
biwords the integers $w_i$ and $w_{i+1}$ must be different. This
implies that the set of tableaux $T$ with $s$ symbols and
associated word $y$ corresponds bijectively to the set of words
$w$ with $1\leq w_1\leq w_2\leq\cdots\leq w_n\leq s$, where the
inequalities are strict in correspondence of the descents of $y$.
Every such word $w$ is uniquely determined by the sequence $\delta
:=w_1-1,w_2-w_1,\ldots,w_n-w_{n-1},s-w_n$, which is a composition
of the integer $s-1$ such that its $i$-th component $\delta_i$ is
at least one whenever $y$ has a descent at the $i$-th position.
For this reason, we can consider the word $\delta'$ defined as
follows:
$$\delta'_i=\left\{\begin{array}{ll}
\delta_i-1 & \textrm{ if $y$ has a descent at the i-th position}\\
\delta_i & \textrm{ otherwise}
\end{array}\right.,$$
which is, of course, a composition of the integer $s-k-1$.\newline

\noindent For example the semistandard tableau on $5$ symbols
$$T=\begin{array}{ccc}1&2&3\\2&3&\\4&&\\5&&\end{array}$$ is associated to the Yamanouchi word $y=1\,1\,2\,1\,2\,3\,4$
with descents at positions $1$ and $3$. In this case, we have:
$$w=1\,2\,2\,3\,3\,4\,5$$
$$\delta=0\,1\,0\,1\,0\,1\,1\,0$$
$$\delta'=0\,0\,0\,0\,0\,1\,1\,0.$$

\noindent We are now in position to prove the following:

\newtheorem{total}[yama]{Theorem}
\begin{total}
\label{abin} The total number of semistandard Young tableaux with
$n$ cells and $k$ symbols is
\begin{equation}a_{n,s}=\sum_{k=0}^{s-1} {n+k\choose k}
y_{n,s-k-1}\label{qualunque}\end{equation} and conversely,
\begin{equation}
y_{n,k}=\sum_{j=1}^{k+1} (-1)^{k-j+1}{n+1\choose
k-j+1}a_{n,j}.\label{qualsiasi}
\end{equation}
\end{total}

\noindent \emph{Proof} The preceding observations show that the
semistandard tableaux with $s$ symbols and associated word $y$ are
in bijection with the compositions of the integer $s-k-1$ into
$n+1$ parts. In other terms, the number of semistandard Young
tableaux with $s$ symbols whose associated reverse Yamanouchi word
$y$ has $k$ descents is
$${n+s-k-1\choose n}.$$ Formula (\ref{qualunque}) follows directly
by these considerations. The second identity can be easily
obtained by inversion.\begin{flushright} $\diamond$
\end{flushright}

\noindent Combining formulae (\ref{a}) and (\ref{qualsiasi}) we
get an explicit expression for the integers $y_{n,k}$:
\newtheorem{expli}[yama]{Corollary}
\begin{expli}
The number $y_{n,k}$ of reverse Yamanouchi words of length $n$
with $k$ descents is
\begin{equation}y_{n,k}=\sum_{j=1}^{k+1} (-1)^{k-j+1}{n+1\choose k-j+1}\sum_{i=0}^{\lfloor\frac{n}{2}\rfloor}{{j \choose 2}+i-1\choose i}
{n+j+2i-1\choose j-1}. \label{explicita}\end{equation}
\end{expli}\begin{flushright} $\diamond$
\end{flushright}

\noindent We remark that Formula (\ref{qualunque}) implies the
following relation between the generating function
$A_n(x)=\sum_{k\geq 0} a_{n,k}x^k$ of the $n$-th row of the matrix
$A$ and the polynomial $I_n(x)$:

\newtheorem{gene}[yama]{Theorem}
\begin{gene}
We have:
$$A_n(x)=\frac{xI_n(x)}{(1-x)^{n+1}}.$$
\end{gene}
\begin{flushright}
$\diamond$
\end{flushright}

\noindent In conclusion of this section, we submit that Theorem
7.19.7 in \cite{stabook} can be rephrased as an immediate
consequence of the described correspondence between reverse
Yamanouchi words and semistandard tableaux. In fact, let $\lambda$
be a partition of the integer $n$ and let $Y(\lambda)$ be the set
of Yamanouchi words whose associated standard tableau has shape
$\lambda$. For every $y\in Y(\lambda)$, we denote by $S(y)$ the
set of semistandard tableaux associated with $y$. The fundamental
quasi-symmetric function $L_y$ can be defined as:

$$L_y(x_1,\ldots, x_m)=\sum_{\scriptsize{\begin{array}{c}1\leq i_1\leq\cdots\leq i_n\leq m\\ i_j<i_{j+1}\ \textrm{if}\ j\in des(y)\end{array}}}x_{i_1}\cdots x_{i_n}.$$
Then, the Schur function $s_{\lambda}(x_1,\ldots,x_m)$ can be
expressed in terms of fundamental quasi-symmetric functions as
follows:
$$s_{\lambda}(x_1,\ldots,x_m)=\sum_{\scriptsize{\begin{array}{c}S\  semistandard\\ sh(S)=\lambda\end{array}}}x^{w(S)}=\sum_{y\in Y(\lambda)}\sum_{S\in S(y)}x^{w(S)}=$$
$$=\sum_{y\in Y(\lambda)}\sum_{\scriptsize{\begin{array}{c}1\leq i_1\leq\cdots\leq i_n\leq m\\ i_j<i_{j+1}\ \textrm{if}\ j\in des(y)\end{array}}}x_{i_1}\cdots x_{i_n}=\sum_{y\in Y(\lambda)}L_y,$$
where $w(S)$ is the content of $S$.

\section{Disproof of the conjecture}

\noindent First of all, we recall a general result appearing in
\cite{uni}:

\newtheorem{rodin}[yama]{Proposition}
\begin{rodin}
\label{pdu} The product $p(x)\cdot q(x)$ of a unimodal polynomial
$p(x)$ and a log-concave polynomial $q(x)$ is unimodal. If $p(x)$
is log-concave, the product $p(x)\cdot q(x)$ is log-concave as
well.
\end{rodin}

\noindent The relation between the sequences $(y_{n,k})$ and
$(a_{n,k})$ described in Theorem \ref{abin} allows us to refute
the log-concavity of the polynomials $I_n(x)$. In fact, we have:

\newtheorem{nolog}[yama]{Theorem}
\begin{nolog}
The polynomials $I_n(x)$ are in general not log-concave.
\end{nolog}

\noindent \emph{Proof} Formula (\ref{qualunque}) shows that the
polynomial
$$p_n(x)=\sum_{k=0}^n a_{n,k}x^k$$ is the product of
the two polynomials $I_n(x)$ and
$$q_n(x)=\sum_{k=0}^n{n+k-1\choose k}x^k.$$ The polynomial
$q_n(x)$ is log-concave. In fact, the condition
$${n+k-2\choose k-1}{n+k\choose k+1}\leq {n+k-1\choose k}$$ is
equivalent to $$\frac{n+k}{k+1}\leq\frac{n+k-1}{k}$$ that holds
for every $n\geq 1$. Hence the log-concavity of $I_n(x)$ would
imply the log-concavity of the polynomial $p_n(x)$, contradicting
Proposition \ref{lapide}.
\begin{flushright}
$\diamond$
\end{flushright}

\noindent In fact, exploiting Formula (\ref{explicita}), we get,
for instance
$$y_{50,1}^2=390625<465570=y_{50,0}\cdot y_{50,2}.$$


\begin{thebibliography}{99}
\bibitem{bbsa} M.Barnabei, F.Bonetti, M.Silimbani, The Eulerian distribution on self evacuated involutions, submitted.

\bibitem{bbsb} M.Barnabei, F.Bonetti, M.Silimbani, The signed Eulerian numbers on involutions, submitted.

\bibitem{bon2} M.Bona, R.Ehrenborg, A combinatorial proof of the log-concavity of the numbers of permutations with $k$ runs,
 \emph{J. Combin. Theory Ser. A} {\bf 90} (2000), no. {\bf 2}, 293--303

\bibitem{comt} L.Comtet, Advanced Combinatorics, Reidel, Dordrecht
(1974).

\bibitem{des} J.D\'esarm\'enien, D.Foata, Fonctions sym\'etriques et
s\'eries hyperg\'eom\'etriques basiques multivari\'ees,
\emph{Bull. Soc. Math. France} {\bf 113} (1985), 3--22.

\bibitem{duk} M.W.B.Dukes, Permutation statistics on involutions, \emph{European J. Combin.} {\bf 28} Issue {\bf 1} (2007), 186-198.

\bibitem{gess} I.M.Gessel, C.Reutenauer, Counting permutations
with a given cycle structure and descent set, \emph{J. Combin.
Theory Ser. A} {\bf 13} (1972), 135-139.

\bibitem{zeng} V.J.Guo, J.Zeng, The Eulerian distribution on involutions is indeed unimodal,
 \emph{J. Combin. Theory Ser. A} {\bf 113} (2006), no. {\bf 6}, 1061--1071.

\bibitem{uni} V.E.Levit, E.Mandrescu, Independence polynomials of well-covered graphs: Generic counterexamples for the unimodality conjecture,
\emph{European J. Combin.} {\bf 27} Issue {\bf 6} (2006), 931-939.

\bibitem{Lot} M.Lothaire, Algebraic combinatorics on words, \emph{Encyclopedia of
Mathematics and its Applications} {\bf 90} Cambridge University
Press, Cambridge (2002).

\bibitem{prod}  H.Prodinger, Some information about the binomial transform, \emph{Fibonacci Quart.} {\bf 32} (1994),
 no. {\bf 5}, 412--415.

\bibitem{schu} I.Schur, $\ddot{\textrm{U}}$ber die rationalen Darstellungen der allgemeinen linearen Gruppe, S'ber, \emph{Akad. Wiss. Berlin} (1927), 58-75, Ges. Abh. III,
68-85.

\bibitem{stabook} R.P.Stanley, Enumerative Combinatorics, Vol. II, \emph{Cambridge
Studies in Advanced Mathematics}, {\bf 62}. Cambridge University
Press, Cambridge (1999).

\bibitem{sta} R.P.Stanley, Log-concave and unimodal sequences in algebra, combinatorics, and geometry, Ann. New York Acad. Sci., {\bf 576} (1989), 500-535.

\bibitem{stem} J.R.Stembridge, Eulerian numbers, tableaux, and the
Betti numbers of a toric variety, \emph{Discrete Math.} {\bf 99}
(1992), 307-320.

\bibitem{stre} V.Strehl, Symmetric Eulerian distributions for
involutions, \emph{S\'eminaire Lotharingien Combinatoire} {\bf 1},
Strasbourg 1980, Publications del l'I.R.M.A. 140/S-02, Strasbourg
1981.

\end{thebibliography}
\end{document}